\documentclass{amsart}
\usepackage{amssymb, mathtools, fullpage, color}
\usepackage[colorinlistoftodos]{todonotes}
\usepackage[colorlinks=true, pdfstartview=FitV,linkcolor=blue,citecolor=blue,urlcolor=blue]{hyperref}

\theoremstyle{definition}
\newtheorem {theorem}{Theorem}[section]
\newtheorem {lemma}[theorem]{Lemma}
\newtheorem {corollary}[theorem]{Corollary}

\newtheorem{definition}[theorem]{Definition}

\newtheorem{example}[theorem]{Example}

\newenvironment{red}{\relax\color{red}}{\relax}
\newenvironment{blue}{\relax\color{blue}}{\hspace*{.5ex}\relax}

\newcommand{\ber}{\begin{red}}
\newcommand{\er}{\end{red}}
\newcommand{\beb}{\begin{blue}}
\newcommand{\eb}{\end{blue}}

\newcommand*\inv[1]{\overleftrightarrow{#1}}

\numberwithin{equation}{section}

\begin{document}

\title[Divisibility in factoradic form]{Divisibility rules for integers presented as permutations}

\date{\today}

\author[T. Oliver]{Thomas Oliver}
\address{University of Westminster, London, U.K.}
\email{T.Oliver@westminster.ac.uk}

\author[A.~Vernitski]{Alexei Vernitski}
\address{University of Essex}
\email{asvern@essex.ac.uk}

\begin{abstract}
In this note, we represent integers in a type of factoradic notation. 
Rather than use the corresponding Lehmer code, we will view integers as permutations.
Given a pair of integers $n$ and $k$, we give a formula for $n\bmod k$ in terms of the factoradic digits, and use this to deduce various divisibility rules.
\end{abstract}


\maketitle

\section{Introduction}\label{s.motivation}

Consider a pair of integers $n$ and $k$.
If $n$ is written in decimal form, then the value of $n\bmod k$ can sometimes be expressed as a linear expression modulo $k$ depending on a fixed subset of the digits of $n$.
For example, $n\bmod 4$ is equal to $2a+b \bmod 4$, where $a$ and $b$ are the two last decimal digits of $n$.
This observation may be used to formulate divisibility rules for $k$.
A list of such rules for the first $1000$ prime numbers is given in \cite{briggs2000simple}.

In this note, we will investigate the analagous statements when $n$ is written in so-called factoradic form. 
Instead of the usual Lehmer code, each non-negative integer will be presented as a permutation on the set $\mathbb{Z}_{\geq0}$ of non-negative integers.
The permutations involved are tame in the sense that they fix all but finitely many integers.
Precise definitions will be given in Section~\ref{s.factoradic}. 

This is, in part, an exercise in data presentation.
When exploring data, the manner in which it is presented changes the nature of the possible experiments, and hence of what can be (machine) learned from them.
Whilst the integers may seem to be a somewhat elementary dataset, there are many questions about them that are considered to be intractable to advanced machine learning architectures.
For example, see the recent case study in deep learning the Mobius function \cite{DLD}.
The authors' experience seems to indicate that it is convenient to represent pure mathematical data as permutation matrices, and then it is especially convenient to represent these permutations using their inversion sets. 
One example of this approach appeared in \cite{lisitsausing}, in which the second named author studied so-called petal diagrams of knots. 
This note presents another example corroborating this heuristic. 
We consider non-negative integers presented as permutations, and we discover that divisibility rules can be easily expressed via the corresponding inversion sets. 
In particular, we will prove the following theorem:
\begin{theorem}\label{t.thm1}[Corollary \ref{c.cor}]
For $n\in\mathbb{Z}_{\geq0}$ and $k\in\mathbb{Z}_{>0}$, we have
$$n \equiv \sum_{i<j} \inv{i, j} \cdot j! \bmod k,$$ 
where the values $\inv{i, j}$ are considered in the inversion set of the $k$-prefix for the factoradic form for $n$.
\end{theorem}
The coefficients $\inv{i, j}$ are essentially the digits of $n$ in factoradic form, and will be explained in section~\ref{s.factoradic}.
Two interesting features of Theorem~\ref{t.thm1} are, firstly, that it provides divisibility formulas depending only on a finite number of factoradic digits, and, secondly, that these formulas are linear.
One can argue that it is not so surprising that the formula depends only on a finite number of digits. 
Indeed, checking divisibility by $2$ for integers in binary form, or checking divisibility by $10$ for integers in decimal form, involves only checking the last digit, and, in a sense, the factoradic form of an integer is $b$-ary for any base $b$.
On the other hand, it is quite surprising that the formulas are linear. 
Although one is used to linear divisibility formulas in decimal arithmetic (see, for example, the afforementioned summary in \cite{briggs2000simple}), the digits $\inv{i, j}$ are of a completely different nature to the usual decimal digits (or any other $b$-ary digits).

We conclude this introduction with a summary of what follows. 
In Section~\ref{s.factoradic}, we define the factoradic form of an integer.
In Section~\ref{s.divisibility} we define the inversion set of an integer in factoradic form, prove Theorem~\ref{t.thm1}, and state some simple divisibility rules. 

\section{What is the factoradic form of a non-negative integer?}\label{s.factoradic}

In this section, we present facts which have been known since the 19th century \cite{laisant1888numeration}.

\begin{definition}
We will refer to a permutation $\sigma$ on $\mathbb{Z}_{\geq0}$ as \textsl{tame} if, for all but finitely many $i$, we have $\sigma(i)=i$.
\end{definition}

In this paper, we will write the tame permutation $\sigma$ as the sequence $(\sigma(0),\sigma(1),\sigma(2),\dots)$.
In particular, the trivial permutation is simply the sequence of non-negative integers in their natural order $(0,1,2,\dots)$.

\begin{example}\label{ex.leadingzeros}
Consider a non-negative integer in its usual decimal notation, for example, $2025$.
This notation should be treated, strictly speaking, as a left-infinite sequence $\dots 0 \dots 02025$, but (almost) all leading zeros can be ignored, as long as all the non-zero digits are included, so the number can be written as $2025$ (or $02025$, or $002025$, and so on).
Similarly, when we work with tame permutations, there is no need to include the entire right-infinite sequence, as long one includes all terms that have been shifted away from their natural position. 
Thus, for example, the tame permutation that swaps $0$ and $1$ and leaves all other numbers in their natural positions can be written as $(1, 0)$, or $(1, 0, 2)$, or $(1, 0, 2, 3)$, etc.
\end{example}

\begin{definition}
The \textsl{factoradic order} on the set of tame permutations is the inverse of the right-to-left lexicographic order. 
\end{definition}

\begin{example}
The non-negative integers in their natural order $(0)$ is the $0$th permutation in factoradic order, and the permutation $(1, 0)$ from Example~\ref{ex.leadingzeros} is the $1$st. 
\end{example}

Note that the factoradic order is linear and, to be more specific, is order-isomorphic to the chain of non-negative integers.

\begin{definition}
The \textsl{factoradic form} of $n\in\mathbb{Z}_{\geq0}$ is the $n$th permutation in the factoradic order of tame permutations. 
\end{definition}

\begin{table}
\centering
\begin{tabular}{|c|c|}
\hline
Decimal & Factoradic\\
\hline
$0$ & $(0, 1, 2, 3)$ \\
$1$ & $(1, 0, 2, 3)$ \\
$2$ & $(0, 2, 1, 3)$ \\
$3$ & $(2, 0, 1, 3)$ \\
$4$ & $(1, 2, 0, 3)$ \\
$5$ & $(2, 1, 0, 3)$ \\
\hline
\end{tabular}
\quad
\begin{tabular}{|c|c|}
\hline
Decimal & Factoradic\\
\hline
$6$ & $(0, 1, 3, 2)$ \\
$7$ & $(1, 0, 3, 2)$ \\
$8$ & $(0, 3, 1, 2)$ \\
$9$ & $(3, 0, 1, 2)$ \\
$10$ & $(1, 3, 0, 2)$ \\
$11$ & $(3, 1, 0, 2)$ \\
\hline
\end{tabular}
\quad
\begin{tabular}{|c|c|}
\hline
Decimal & Factoradic\\
\hline
$12$ & $(0, 2, 3, 1)$ \\
$13$ & $(2, 0, 3, 1)$ \\
$14$ & $(0, 3, 2, 1)$ \\
$15$ & $(3, 0, 2, 1)$ \\
$16$ & $(2, 3, 0, 1)$ \\
$17$ & $(3, 2, 0, 1)$ \\
\hline
\end{tabular}
\quad
\begin{tabular}{|c|c|}
\hline
Decimal & Factoradic\\
\hline
$18$ & $(1, 2, 3, 0)$ \\
$19$ & $(2, 1, 3, 0)$ \\
$20$ & $(1, 3, 2, 0)$ \\
$21$ & $(3, 1, 2, 0)$ \\
$22$ & $(2, 3, 1, 0)$ \\
$23$ & $(3, 2, 1, 0)$ \\    
\hline
\end{tabular}
\\~\\
\caption{Factoradic forms of $0, 1, \dots, 23$.}
\label{tab:factoradic-examples}
\end{table}

\begin{example}
The factoradic forms of the nunbers in the set $\{0,1,\dots,23\}$ are presented in Table \ref{tab:factoradic-examples}.
The number $23$ is not chosen accidentally. 
Indeed, the set of non-negative integers in factoradic form possesses a certain periodicity, with the periods being factorials, and $23 = 4!-1$. 
The meaning of this will be clarified in Lemma~\ref{lemma}.
\end{example}

\begin{example}\label{ex.linkwithfactorials}
The factoradic form of $n$ is connected to a (unique) expression for $n$ of the form
\begin{equation}\label{eq.factorialsum}
n = \sum_{i\geq 0} a_i \cdot i!,
\end{equation}
where each coefficient $a_i$ is an integer in the set $\{0, 1, \dots, i\}$. 
For example, consider the tame permutation $(2, 3, 0, 1)$, which represents the decimal number $16$. 
Starting from the right, we consider the positions in which the entry is smaller than the index. 
In the natural order, we would expect $3$ to be in position $3$, but we have $1$ there. 
This indicates that the term $(3-1)\cdot 3!=2 \cdot 3!$ will appear in equation~\eqref{eq.factorialsum}. 
Now let us proceed to entries $2, 3, 0$ in positions $0, 1, 2$. 
If these entries were in the increasing order (approximating the natural order), they would have been in the order $0, 2, 3$. 
Instead of this, in position $2$ we have not $3$, but $0$, that is, the second smallest number of those that remain. 
This indicates that a term $(2-0)\cdot 2!=2\cdot 2!$ will appear in equation~\eqref{eq.factorialsum}. 
Proceeding, we are left with $2, 3$ in positions $0, 1$, as per their natural order.
Adding the non-trivial terms together, we obtain $2 \cdot 3! + 2 \cdot 2! = 16$. 
\end{example}

\begin{definition}
For $s\geq0$, the $s$-prefix of a tame permutation is the first $s$ elements in the sequence $(\sigma(0),\sigma(1),\sigma(2),\dots)$, that is, the non-negative integers in positions $0, 1, \dots, s-1$. 
\end{definition}

For example, Table \ref{tab:factoradic-examples} lists $4$-prefixes of all those tame permutations that fix all integers $\geq4$. 
From the discussion above, it is obvious that all non-negative integers up to $s!-1$, and only they, are presented by tame permutations that fix all integers $\geq s$.

\section{Using inversions to describe divisibility}\label{s.divisibility}

\begin{definition}
Consider an integer $s\geq1$ and let $\sigma$ be a permutation on a set $\{a_i:i=0,\dots,s-1\}$ non-negative integers such that $a_0 < \cdots < a_{s-1}$. 
For $i<j\in\{0,\dots,s-1\}$, we say that $\sigma$ is an \textit{inversion} of $a_i$ and $a_j$ if $a_i=\sigma(a_{i'})$ and $a_j=\sigma(a_{j'})$ for some $j' < i'$. 
We introduce the value
$$\inv{i, j}=\begin{cases}1,&\sigma \text{ is inversion of }a_i\text{ and } a_j,\\0,&\text{otherwise}.\end{cases}$$ 
We will refer to the values $\inv{i, j}$, for all $i < j$, as the \textit{inversion set} of $\sigma$. 
\end{definition}
Note that, strictly speaking, the inversion set is not a `set' but a function from the set $\{(i, j):i<j\}$ to $\{0, 1\}$. 
We also caution that the inversion set has nothing to do with the inverse of the permutation matrix.

\begin{lemma}\label{lemma}
Consider $m, n, s\in\mathbb{Z}_{\geq0}$. 
If $m \equiv n \bmod s!$, then the $s$-prefixes for the factoradic forms of $m$ and $n$ have the same inversion sets.
\end{lemma}

Lemma~\ref{lemma} expresses what we saw in the examples in Table \ref{tab:factoradic-examples}. 
Although the numbers standing in the first $s=2$ or $s=3$ positions might change in a complicated way from one part of Table \ref{tab:factoradic-examples} to another, one can spot that the inversion set of the $s$-prefixes repeats regularly, namely, with period $s!$. For example, one can easily check that the inversion sets of $3$-prefixes of the tame permutations representing numbers $0$, $6$, $12$ and $18$ coincide.

\begin{theorem}\label{t.thm}
Consider $n,s\in\mathbb{Z}_{\geq0}$.
If $n < s!$, then 
$$n = \sum_{i<j} \inv{i, j} \cdot j!,$$ 
where the values $\inv{i, j}$ are considered in the inversion set of the $s$-prefix of the factoradic form for $n$.
\end{theorem}

Theorem~\ref{t.thm} may be proved by arguing as in Example~\ref{ex.linkwithfactorials}. 
Theorem~\ref{t.thm} implies the following Corollary, previously stated as Theorem~\ref{t.thm1}, which expresses whether or a non-negative integer $n$, represented in the factoradic form, is divisible by another integer $k$.

\begin{corollary}\label{c.cor}
For $n\in\mathbb{Z}_{\geq0}$ and $k\in\mathbb{Z}_{>0}$, we have
$$n \equiv \sum_{i<j} \inv{i, j} \cdot j! \bmod k,$$ 
where the values $\inv{i, j}$ are considered in the inversion set of the $k$-prefix for the factoradic form for $n$.
\end{corollary}

We conclude this section with some simple examples of how Corollary~\ref{c.cor} can be applied.
\begin{example}
We apply Corollary~\ref{c.cor} with the following values of $k$:
\begin{enumerate}
\item When $k=2$, we deduce that $n\in\mathbb{Z}_{\geq0}$ is divisible by $2$ if and only if $\inv{0, 1} = 0.$
In other words, an integer in factoradic form is even if the first two terms in the sequence are in their natural order.
Looking at Table \ref{tab:factoradic-examples}, one may have conjectured that this is true. 
\item When $k=3$, we deduce that $n\in\mathbb{Z}_{\geq0}$ is divisible by $3$ if and only if the following is divisble by $3$:  
$$\inv{0, 1} - \inv{0, 2} - \inv{1, 2}.$$ 
It is easy to verify that this is true on examples of the numbers in Table \ref{tab:factoradic-examples}, but you are unlikely to conjecture this easily by considering Table \ref{tab:factoradic-examples}.
\item When $k=4$, we deduce that $n\in\mathbb{Z}_{\geq0}$ is divisible by $4$ if and only if the following is divisble by $4$: 
$$\inv{0, 1} + 2\left(\inv{0, 2} + \inv{0, 3} + \inv{1, 2} + \inv{1, 3} +  \inv{2, 3}\right).$$ 
\item When $k=5$, we deduce that $n\in\mathbb{Z}_{\geq0}$ is divisible by $5$ if and only if the following is divisble by $5$: 
$$\inv{0, 1} + 2\left(\inv{0, 2}\right) + \inv{0, 3} - \inv{0, 4} + 2\left(\inv{1, 2}\right) + \inv{1, 3}  - \inv{1, 4} +  \inv{2, 3} -  \inv{2, 4} -  \inv{3,4}.$$ 
\item When $k=6$, we deduce that $n\in\mathbb{Z}_{\geq0}$ is divisible by $6$ if and only if the following is divisble by $6$: 
$$\inv{0, 1} + 2(\inv{0, 2}+\inv{1,2}).$$ 
\end{enumerate}
\end{example}

\end{document}